\documentclass[12pt,twoside]{amsart}\usepackage{amssymb}
\usepackage{enumerate}

\usepackage{amsmath}
\usepackage[all]{xy}

\usepackage[mathscr]{eucal}
\usepackage{mathrsfs}

\newcommand{\cali}[1]{\mathscr{#1}}
\newcommand{\Hc}{{\cali H}}
\newtheorem{thmspec}{\relax}
\newtheorem{theorem}{Theorem}[section]

\newtheorem{lem}[theorem]{Lemma}

\theoremstyle{definition}

\theoremstyle{remark}

\numberwithin{equation}{section} \tolerance = 10000
\def \Bbb{\mathbb}

\def\onto{{\kern3pt\to\kern-8pt\to\kern3pt}}

\def\<{\langle}
\def\>{\rangle}
\def\|{{\ |\ }}

\def\onto{\twoheadrightarrow}

\def\-{\underline}

\def\R{\Bbb R}
\def\C{\Bbb C}

\def\X{\Bbb X}

\def\<{\langle}
\def\>{\rangle}

\catcode`\@=11
\def\serieslogo@{\relax}
\def\@setcopyright{\relax}
\catcode`\@=12

\title[Envelope of holomorphy for boundary cross sets]
{ Envelope of holomorphy for boundary cross sets
 }

\begin{document}

\author{Peter Pflug}
\address{Peter Pflug\\
Carl von  Ossietzky Universit\"{a}t Oldenburg \\
Fachbereich  Mathematik\\
Postfach 2503, D--26111\\
 Oldenburg, Germany}
\email{pflug@mathematik.uni-oldenburg.de}

\author{Vi{\^e}t-Anh  Nguy\^en}
\address{Vi{\^e}t-Anh  Nguy\^en\\
Mathematics Section\\
The Abdus Salam international centre
 for theoretical physics\\
Strada costiera, 11\\
34014 Trieste, Italy} \email{vnguyen0@ictp.trieste.it}

\date{}


\begin{abstract}
 Let $D\subset \C^n,$ $G\subset \C^m$ be   open sets, let
  $A$ (resp. $B$) be  a subset of the boundary $\partial D$ (resp.
  $\partial G$) and let $W$ be the  $2$-fold boundary cross $((D\cup A)\times B)\cup (A\times(B\cup G)).$
   An open  subset $X\subset\C^{n+m} $ is said to be the ``envelope of holomorphy"
  of $W$ if it  is, in some sense, the maximal open  set with the following property:  Any function
  locally bounded  on $W$ and separately
  holomorphic
  on $(A\times G) \cup (D\times B)$ ``extends" to a holomorphic  function    defined on $X$
   which admits the boundary values $f$   a.e. on  $W.$ In this work
   we will determine the envelope of holomorphy
  of some boundary crosses.
     \end{abstract}
\maketitle

\section{ Introduction}

In a series of articles \cite{pn1,pn2,pn3} the authors establish
various ``boundary cross theorems". These results      deal with
the continuation of holomorphic functions of several complex
variables which are defined on some boundary crosses. The first
theorem of this type was discovered and proved by
Malgrange--Zerner \cite{ze}.

However,  the question naturally arises whether all these theorems
are optimal. More precisely,  are the extension domains  in these
theorems   always maximal? In other words, are they always
``envelopes of holomorphy"? In this work we investigate this
question. We will show that under some conditions our boundary
cross theorems are optimal.

\subsection{Plurisubharmonic measures}
Let $\Omega\subset \C^n$ be an open set.  For any function $u$
defined on $\Omega,$ let
\begin{equation*}
 \hat{u}(z):=
\begin{cases}
u(z),
  & z\in   \Omega,\\
 \limsup\limits_{\Omega\ni w\to z}u(w), & z \in \partial \Omega.
\end{cases}
\end{equation*}
For a set  $A\subset \overline{\Omega}$ put
\begin{equation*}
h_{A,\Omega}:=\sup\left\lbrace u\ :\  u\in\mathcal{PSH}(\Omega),\
u\leq 1\ \text{on}\ \Omega,\
   \hat{u}\leq 0\ \text{on}\ A    \right\rbrace,
\end{equation*}
where $\mathcal{PSH}(\Omega)$ denotes the cone   of all functions
plurisubharmonic on $\Omega.$

The {\it plurisubharmonic measure} of $A$ relative to $\Omega$ is
the function  $\omega(\cdot,A,\Omega)\in \mathcal{PSH}(\Omega)$
defined by
\begin{equation*}
\omega(z,A,\Omega):=h^{\ast}_{A,\Omega}(z),\qquad z\in\Omega,
\end{equation*}
where  $h^{\ast}$ denotes the upper semicontinuous regularization
of the function $h.$

If $n=1$ and $A\subset\partial \Omega,$ then
$\omega(\cdot,A,\Omega)$ is also called
 the {\it  harmonic measure} of $A$ relative to $\Omega.$
 In this case,  $\omega(\cdot,A,\Omega)$
is a harmonic function.
\subsection{Cross, separate holomorphicity.}

For open sets $D\subset \C^n,$  $G\subset \C^m$ and  subsets
$\varnothing\not = A\subset \overline{D},$  $\varnothing\not =
B\subset \overline{G},$   we define
 the {\it cross} $W$ and {\it its interior} $W^{\text{o}}$  as
\begin{eqnarray*}
W &=&\X(A,B; D,G)
:= ((D\cup A)\times B)\cup(A\times (G\cup B)),\\
W^{\text{o}} &=&\X^{\text{o}}(A,B;D,G) :=  (D \times
B)\cup(A\times G ).
\end{eqnarray*}
If $A\subset \partial D$ and $B\subset\partial G$  (resp.
$A\subset  D$ and $B\subset G$), then $W$ is called a {\it
boundary cross}  (resp.  a {\it classical cross}).

 For a cross $W:=\X(A,B;D,G)$
 let
 \begin{equation*}
 \widehat{W}^{\text{o}} =\widehat{\X}^{\text{o}}(A,B;D,G) :=\left\lbrace  (z,w)\in D \times
 G :\ \omega(z,A,D)+\omega(w,B,G)  <1
\right\rbrace
\end{equation*}
and
\begin{equation*}
\widehat{W}=\widehat{\X} (A,B;D,G):= W\cup  \widehat{W}^{\text{o}}
.
\end{equation*}

We say that a function $f:\ W\longrightarrow \C$ is  {\it
separately holomorphic} on $W^{\text{o}} $ and write
$f\in\mathcal{O}_s(W^{\text{o}}),$   if
 for any  $a\in A$
 the function $f(a ,\cdot )|_{G}$ is holomorphic  on $G,$  and for any $b\in B$
   the function $f( \cdot,b )|_{D}$ is holomorphic  on $D.$
%In the sequel, $\mathcal{C}(W)$ denotes the space of all complex-valued continuous functions on $W.$

We say that a function $f:\ W\longrightarrow \C$ is  {\it
separately continuous} and write $f\in\mathcal{C}_s(W),$   if
 for any  $a\in A$   and for any $b\in B,$
 the functions $f(a ,\cdot )$ and $f( \cdot,b )$ are continuous.

For an open set $\Omega\subset\C^n,$ $\mathcal{O}(\Omega)$ denotes
the space of all holomorphic functions
 on $\Omega.$ % For a function $f:\ M\longrightarrow\C$ and  a subset $N\subset M,$  $f|_N$ denotes the
 %the function obtained by restricting $f$ on $N.$

%
%
\subsection{Motivations for our work and envelope of holomorphy of a boundary cross}

We like to formulate   the {\it  boundary cross theorems}  in one
and higher dimensional contexts (see \cite{pn1,pn2,pn3}).

%For the one dimensional case we refer the reader to  the notation and terminology
% in Subsections 2.1--2.3 in the work \cite{pn2}.
  A {\it (Jordan) curve} in $\C$
  is the image $\mathcal{C}:=\{\gamma(t):\ t\in[0,1]\}$ of a continuous one-to-one map    $
\gamma:\ [0,1]\longrightarrow  \C.$ The    interior of the curve
$\mathcal{C}$ given by $ \{\gamma(t):\ t\in(0,1)\}$ is said to be
an {\it open (Jordan) curve}.
   A {\it Jordan domain}  is the image $\{\Gamma(t), \ t\in E\}$
of  a   one-to-one continuous  map    $ \Gamma:\
\overline{E}\longrightarrow  \C,$ where, in this work, $E$ denotes
the open  unit disc in $\C.$
 A {\it closed (Jordan) curve} is the boundary  of a Jordan domain.
 An open set   $D\subset  \C$   is said to be {\it Jordan-curve-like  at a point}
$\zeta\in\partial D$ if there is a  Jordan domain $U$ such that
$\zeta\in U $  and  %  $U\setminus \partial D=V_1\cup V_2$ for two
%Jordan domains $V_1,$ $V_2$ and $\overline{V_1}\cap\overline{
%V_2}\cap U=
 $U\cap\partial D$ is an open (Jordan) curve.

 Let $D\subset \C,$ $ G\subset \C$ be two open sets and
  $A$ (resp. $B$)  a subset of  $\partial D$ (resp.
  $\partial G$) such that
      $D$ (resp. $G$) is   Jordan-curve-like
   at every point of  $A$ (resp. $B$), and  let $f:\ W\longrightarrow \C$ be a function.
 We can define as in Subsections 2.1--2.3 of  \cite{pn2} various
notions and terminology:  Jordan-measurable sets, sets of positive
length, sets of zero length, Jordan-measurable functions,  the
angular limit, the set of all locally regular points $A^{\ast}$
(resp. $B^{\ast}$)
 relative to $A$  (resp.  $B$),  almost everywhere (a.e.)  etc.

 Theorem A in \cite{pn2} may be restated, in a simple form, %,  in a simple form (i.e. $X=Y=\C$),
  as follows:
    %Section 2 in \cite{pn3}. For example, $\mes$ denotesthe linear measure.

\renewcommand{\thethmspec}{Theorem 1}
\begin{thmspec}
 We keep the hypotheses and notation of the previous paragraph.
  Suppose in addition that $A$ and $B$ are  of positive  length
  and that $f$ verifies the following properties:
%Let  $f:\ W\longrightarrow \C$ be   such that:
\begin{itemize}
\item[ (i)] $f$ is locally bounded on $W$ and  $f\in
\mathcal{O}_s(W^{\text{o}});$ \item[(ii)] $f|_{A\times B}$ is
Jordan-measurable; \item[ (iii)]   for any $a\in  A$  (resp. $b\in
B$),
 the holomorphic function $f(a,\cdot)|_{G}$  (resp.  $f(\cdot,b)|_{D}$) has the angular  limit
  $f_1(a,b)$ at $b$ for a.e. $b\in B$  (resp.  $f_2(a,b)$ at  $a$ for a.e. $a\in A$)
and   $f_1=f_2=f$  a.e. on $ A\times B.$
\end{itemize}

Then  there exists a unique function
$\hat{f}\in\mathcal{O}(\widehat{W}^{\text{o}})$ with the following
property:   There are  subsets $\tilde{A}\subset A\cap A^{\ast}$
and $ \tilde{B}\subset B\cap B^{\ast}$  such that
\begin{itemize}
\item[1a)]
 the sets  $A\setminus  \tilde{A}$ and    $B\setminus  \tilde{B}$
are of zero length;%\footnote{ Under this condition it follows from   Part 1) of
%Theorem 4.6 in \cite{pn2}  that  $\tilde{A}\subset\tilde{A}^{\ast}$  and  $\tilde{B}\subset\tilde{B}^{\ast}$   .}
 \item[1b)] $\hat{f}$  admits the angular limit $f(\zeta,\eta)$ at every point
  $(\zeta,\eta)\in(\tilde{A}\times G)\bigcup( D\times\tilde{B}).$
\end{itemize}
\end{thmspec}
In fact, this theorem  was formulated  in \cite{pn2} in a more
general context:
 $D$ and $G$ are open sets of arbitrary complex manifolds of dimension $1$ countable at infinity.

For the higher dimensional case we recall the following
terminology
 from Section 2 in  \cite{pn3}.
 Let  $D\subset\C^{n}$    be a nonempty open set, and
  $A$   a nonempty relatively open subset of  $\partial D.$
  Then $A$ is said to be a {\it  topological hypersurface (in $\C^n\equiv\R^{2n}$)}  if,
  for every $a\in A$ there exist an open neighborhood $V$ of $a,$
an open subset $U\subset \R^{2n-1},$  a continuous function $h:\
U\longrightarrow\R$ and an integer $j:$ $1\leq j\leq 2n$   such
that
\begin{multline*}
V\cap A=\left\lbrace z=(x_1,\ldots,x_{2n})\in \R^{2n}:\  x_j=h(x_1,\ldots,x_{j-1},x_{j+1},\ldots,x_{2n}),\right.\\
\left.  (x_1,\ldots,x_{j-1},x_{j+1},\ldots,x_{2n})\in U
     \right\rbrace.
\end{multline*}

 The Main Theorem in \cite{pn1,pn3} may be restated,  in a simple form, % (i.e. $X=\C^n$ and $Y=\C^m$),
  as follows:

\renewcommand{\thethmspec}{Theorem 2}
  \begin{thmspec}
Let  $D\subset\C^{n},$   $G\subset\C^m$  be two nonempty open
sets, let
  $A$ (resp. $B$) be a nonempty relatively open subset of  $\partial D$ (resp.
  $\partial G$).  Suppose in addition that
  $A$ and $B$ are topological hypersurfaces.
 Let  $f:\ W\longrightarrow \C$ be   such that:
\begin{itemize}
\item[ (i)]  $f\in\mathcal{C}_s(W)\cap
\mathcal{O}_s(W^{\text{o}});$ \item[(ii)]   $f$ is locally bounded
on $W;$ \item[ (iii)]   $f|_{A\times B}$ is continuous.
\end{itemize}

Then  there exists a unique function $\hat{f}\in
\mathcal{O}(\widehat{W}^{\text{o}})$ such that
 \begin{equation*}
 \lim\limits_{(z,w)\to(\zeta,\eta),\ (z,w)\in  \widehat{W}^{\text{o}}}
\hat{f} (z,w) =f(\zeta,\eta),\qquad (\zeta,\eta)\in W.
 \end{equation*}
 \end{thmspec}

In fact, this theorem  was formulated  in \cite{pn3} in its full
generality:
 $D$ and $G$ are open sets of arbitrary complex manifolds.

These results lead to the following concept.

\renewcommand{\thethmspec}{Definition 1}
\begin{thmspec}
 Let  $D,\ G,\ A, \  B $  be as  in the hypothesis of    Theorem 1 (resp. Theorem 2)
 and $W:=\X(A,B;D,G).$
We say that   $\widehat{W}^{\text{o}} $ is the {\bf envelope of
holomorphy} of the boundary cross $W$ if there do not exist
nonempty open sets $U_1,\ U_2 \subset \C^2$  (resp. $\C^{n+m}$)
with $U_2$ connected, $U_2\not\subset \widehat{W}^{\text{o}},$
$U_1\subset U_2\cap \widehat{W}^{\text{o}},$ such that for every
$f:\ W\longrightarrow\C$  which satisfies (i)--(iii) of Theorem 1
(resp. Theorem 2),
 there is a function $h\in\mathcal{O}(U_2)$
such that  $h=\hat{f}$  on $U_1,$  where  $
\hat{f}\in\mathcal{O}(\widehat{W}^{\text{o}})$ is the unique
function given by Theorem 1 (resp. Theorem 2).
 \end{thmspec}
 %Finally, throughout the paper,
%the notation  $\vert f\vert_M$ denotes $\sup_M \vert f\vert.$

The purpose of this article is to investigate the question whether
$\widehat{W}^{\text{o}}$   in  Theorem 1 and 2
 is always the envelope of holomorphy.
This  problem is motivated by the  work of Alehyane--Zeriahi
\cite{az},  where  the envelope of holomorphy of a
 classical cross (i.e. $A\subset D, B\subset G$),
 $D,\ G$ are
 subdomains of Stein manifolds,  has been identified.
 See also  \cite{nv} for further generalizations.

\smallskip

\indent{\it{\bf Acknowledgment.}}   The second author wishes to
express his gratitude to the Max-Planck Institut f\"{u}r
Mathematik in Bonn (Germany) for its hospitality and its support.
The  research was partially supported by DFG grant no. 227/8-2.

\section{Statement of the   results}
%
%
%
%First recall from Subsection 2.1 in \cite{pn2} that given an open set $D\subset \C$
%and a point $\zeta\in\partial D,$    then $\zeta$ is said to be {\it   of type 1}
%if there is an open  neighborhood $V$ of $\zeta$ such that $V\cap D$ is a Jordan
%domain.
%In order to state the first main result we need to introduce some
%more notation and terminology.

 Let $D\subset \C$   be an open set which is
Jordan-curve-like at a point $\zeta\in\partial D.$ Then  $\zeta$
is said to be {\it of type 1} if there is a neighborhood $V$ of
$\zeta$ such that $V\cap D$ is a Jordan domain.  Otherwise,
$\zeta$ is said to be {\it   of type 2}. We easily see  that if
$\zeta$ is of type 2, then there are an open neighborhood $V$ of
$\zeta$ and two Jordan domains $V_1,$ $V_2$  such that  $V\cap
D=V_1\cup V_2.$
 A (Jordan) curve or an open (Jordan) curve or a closed curve  $\mathcal{C}\subset \partial D$ is said to be {\it of type 1} (resp. {\it type 2})
 if all points of $\mathcal{C}$ are of type 1 (resp. type 2).

The following  simple example (see   Subsection 2.1 in \cite{pn2})
may clarify the above definitions.

\smallskip

\noindent {\bf  Example 1.} Let $H$ be the open square in $\C$
whose four vertices are $1+i,$ $-1+i,$ $-1-i,$ and $1-i.$ Define
the domain
\begin{equation*}
D:=H\setminus \left [-\frac{1}{2},\frac{1}{2}\right].
\end{equation*}
Then   $D$ is Jordan-curve-like on $\partial H\cup  \left
(-\frac{1}{2},\frac{1}{2}\right).$ Every point of $\partial H$ is
of type 1 and every point of
$\left(-\frac{1}{2},\frac{1}{2}\right)$ is of type 2. In other
words, $\partial G$ is a closed curve of type 1 and
$\left(-\frac{1}{2},\frac{1}{2}\right)$ is an open curve of type
2.

 We continue with another example  showing that in Theorem 1
 above
  $\widehat{W}^{\text{o}}$   is, in general,  not the envelope of holomorphy of
$W.$

\smallskip

\noindent {\bf  Example 2.} Let $D$ be as in Example 1, let
$A:=\left(-\frac{1}{4},\frac{1}{4}\right),$ and $G:=E,$
$B:=\partial G.$ Then a direct computation shows that
\begin{equation*}
\widehat{W}^{\text{o}} =\widehat{\X}^{\text{o}}(A,B;D,G) =D\times
G.
\end{equation*}
Let $f$ be an arbitrary function satisfying the hypothesis of
Theorem 1, and let $\hat{f}\in
\mathcal{O}(\widehat{W}^{\text{o}})$  and $\widetilde{A}\subset A$
be as in the conclusion this theorem. Therefore, using the
Lindel\"{o}f Theorem, we have
\begin{equation*}
\lim\limits_{y\to 0^{+}}\hat{f}(x+iy)=f(x)=\lim\limits_{y\to
0^{-}}\hat{f}(x+iy),\qquad    x\in\widetilde{A}.
\end{equation*}
 On the other hand, using the fact that the one dimensional Hausdorff measure of
$A\setminus \widetilde{A}$ is zero, the hypothesis that $f$ is
locally bounded on $W$ and applying Two-Constant Theorem, one can
prove that for every $(\zeta,\eta)\in A\times G,$ there are open
neighborhoods $U$ of $\zeta$ in $D\cup A$ and $V$ of $\eta$ in $G$
such that
\begin{equation*}
\sup\limits_{(U\setminus A)\times V}\vert \hat{f}(z,w)\vert
<\infty.
\end{equation*}
Consequently, by Morera's Theorem, $\hat{f}$ extends
holomorphically through all points of $A\times G.$ Therefore, one
can take $U_1:=D\times G$ and $U_2:=(D\cup A)\times G$ in
Definition 1. Hence, $\widehat{W}^{\text{o}}$   is not the
envelope of holomorphy of $W.$

This discussion leads us to the following

 \renewcommand{\thethmspec}{Definition 2}
  \begin{thmspec}
  Let $D\subset \C$  be an open set and $A\subset\partial D$ be such that $D$ is Jordan-curve-like at
  every point of $A.$
  A point  $a\in A$ is said to be {\bf an extendible point} of $A$
   if there exists an open  (Jordan) curve $ \mathcal{C}:=\{\gamma(t):\ t\in  (0,1)
   \}\subset\partial D$
   such that
   \begin{itemize}
   \item[$\bullet$]
    $a\in\mathcal{C},$
  \item[$\bullet$] the open curve
 $ \mathcal{C}$  is of type 2,
\item[$\bullet$]
  $\mathcal{C}\setminus A$ is of zero length.
   \end{itemize}
   \end{thmspec}

Now we are able to state the first result.
\renewcommand{\thethmspec}{Theorem A}
  \begin{thmspec}  Let  $D,\ G,\ A, \  B $  be as  in the hypothesis of    Theorem 1
 and $W:=\X(A,B;D,G).$   Suppose in addition that $D$  (resp. $G$) is Jordan-curve-like at all points of $\overline{A}$
  (resp. $\overline{B}$)  and that  $A$ and $B$ do not possess any extendible points.
 Then  $\widehat{W}^{\text{o}}$   is   the envelope of holomorphy of $W.$
  \end{thmspec}

  For the higher dimensional situation we need to introduce a new  concept.
 % Recall from   \cite{wa}  that
 % a domain  $D\subset \C^n$ is called {\it strongly regular at} $z\in\partial D$  if % for every number $t>0$
 % there exists a neighborhood $U$ of $z$
 % and  a negative function $u\in \mathcal{PSH}(D\cap U)$ such that
 %  $\lim\limits_{w\to z,\ w\in D} u(w)=0$   and  for every number $t>0,$
 %  there are $0<\tau<t$ and $\epsilon>0$  such that
 % such that  $u<-\epsilon$ on $ (D\cap U)\setminus \B(z,\tau),$  where
 % $\B(z,\tau)$ is the ball with center $z$ and radius $\tau.$
 % A domain is called {\it strongly regular} if it is strongly regular at every point $z\in\partial D.$

 % Observe that if either $D $ is strongly pseudoconvex or smooth pseudoconvex of finite type, then it is strongly regular.

 \renewcommand{\thethmspec}{Definition 3}
  \begin{thmspec}
  A pair $(A,D)$ of an open set $D\subset \C^n$ and a relatively  open set $A$ of $\partial D$ is said to satisfy
  hypothesis ($\Hc$) if there exists a sequence of pseudoconvex open sets $(D_k)_{k=1}^{\infty}\subset
   \C^n$ such that
  $D_k\cap \overline{D}=D\cup A$ and   $\bigcap\limits_{k=1}^{\infty} D_k= D\cup A.$
 % there exists a decreasing  sequence of pseudoconvex open sets $(D_k)_{k=1}^{\infty}\subset \C^n$  such that
 %  $D_k\cap\overline{D}=D\cup A,$ $k\geq 1,$ and   $\bigcap\limits_{k=1}^{\infty} D_k= D\cup A.$
   \end{thmspec}

  % it fulfills the following conditions:
  %\begin{itemize}
  %\item[(i)]  $D$ is strongly regular;
  %\item[(ii)] $A$ can be written in the form of a finite disjoint unions
  %$A=\bigcup\limits_{j=1}^N A_j,$ where  $A_j$ are connected  relatively open subsets  $\partial D,$  $j=1,\ldots,N$
  %and $\partial A$
  % (with respect to $\partial D$) is pluripolar (in $\C^n$);
  %\item[(iii)]  for every open neighborhood $U$ of $\overline{A}$ there is a pseudoconvex
  %open set $\Omega$ such that   $D\cup\overline{A}\subset \Omega\subset D\cup U.$
  %\end{itemize}

  Now  we are ready to state the second result.
  \renewcommand{\thethmspec}{Theorem B}
  \begin{thmspec}
  Let  $D,\ G,\ A, \  B $  be as  in the hypothesis of     Theorem 2
 and $W:=\X(A,B;D,G).$ Suppose in addition that
 the pairs $(A,D)$ and $(B,G)$  satisfy hypothesis ($\Hc$).  Then $\widehat{W}^{\text{o}}$   is
   the envelope of holomorphy of $W.$
  \end{thmspec}

  % It is  worthy to   say some words about Definition 2.
   %Since  $W\subset \partial \widehat{W}^{\text{o}},$
%    one can check, using Definition 1, that $ \widehat{W}^{\text{o}}$ is the envelope of holomorphy only if
 %  $D$  and $G$  are pseudoconvex.
   Here  is a simple  sufficient condition. % for  hypothesis (H).

  \renewcommand{\thethmspec}{Proposition C}
  \begin{thmspec}
  Let $D\subset \C^n$ be a pseudoconvex domain  and  $A\subset \partial D$ a relatively  open subset.
  Suppose in addition that  $A$ can be written in the form
 $A=\bigcup\limits_{j\in J} A_j,$ where
 \begin{itemize}
 \item[(i)]  $A_j$ is a connected  relatively open and bounded subset of $\partial D$
 for $j\in J;$
\item[(ii)]  $\overline{A_j}\cap \overline{\bigcup\limits_{k\in
J\setminus\{j\}}A_k}=\varnothing$ for  $j\in J;$ \item[(iii)]
$\bigcup\limits_{j\in J} \overline{A_j}$  is contained in the set
of
 $\mathcal{C}^2$ smooth, strongly pseudoconvex points
of $\partial D.$
% \item[(iv)]    $\partial A_j$ is a
%$\mathcal{C}^2$ smooth $(2n-2)$ real dimensional manifold in
%$\partial D,$
 % $j=1,\ldots,N.$
\end{itemize}
 Then the pair $(A,D)$   satisfies hypothesis ($\Hc$).
 \end{thmspec}

  It seems to be of interest to find weaker conditions than hypothesis ($\Hc$)  for  Theorem B  to be true.  One may also seek to determine the envelope of holomorphy
  of boundary cross sets where some singularities are allowed or in the manifold context. The  problem
  of determining the envelope of holomorphy of classical cross sets  with singularities has been
  studied in many works (see \cite{jp2,jp3,jp4,nv} and the references therein).

\section{The proofs}
The main idea of the proofs  is contained in the following  lemma.
  \begin{lem}\label{lem1}
   Let  $D,\ G,\ A, \  B $  be as  in the hypothesis of    Theorem A (resp. Theorem B)
 and $W:=\X(A,B;D,G).$
  Assume that for every nonempty open connected  set $U\subset \C^2$  (resp. $\C^{n+m}$) such that
$U\not\subset \widehat{W}^{\text{o}}$  and $U\cap
\widehat{W}^{\text{o}}\not=\varnothing.$
   there is a pseudoconvex open set $\Omega$ in $\C^2$ (resp. $\C^{n+m}$) such that
 $\widehat{W}\subset\Omega$ and $      U\cap\partial\Omega\not=\varnothing.$ Then
 $\widehat{W}^{\text{o}}$ is the envelope of holomorphy of $W.$
  \end{lem}
  \begin{proof} Assume that $\widehat{W}^{\text{o}}$ is not the envelope of holomorphy of $W.$
   Then there are non empty
open connected sets $U_1\subset U_2$,
$U_1\subset\widehat{W}^{\text{o}}$ but
$U_2\not\subset\widehat{W}^{\text{o}}$  such that for every  $f:\
W\longrightarrow\C$  which satisfies (i)--(iii) of Theorem 1
(resp. Theorem 2),
 there is a function $h\in\mathcal{O}(U_2)$
such that  $h=\hat{f}$  on $U_1,$  where  $
\hat{f}\in\mathcal{O}(\widehat{W}^{\text{o}})$
 is the unique function given by Theorem 1 (resp. Theorem 2).

% Fix points $w_1\in V_1$ and $w_2\in
%V_2\setminus\widehat{W}^{\text{o}}$ and connect them with a curve
%$\gamma:[0,1]\to U_2.$  Observe that $\dist(\
%\gamma([0,1]),\partial U_2)>\epsilon>0$. Define $U_2:=\cup_{0\leq
%t\leq 1}\B(\gamma(t),\frac{\epsilon}{2})$ and $U_1:=U_2\cap V_1$. Then
%the pair $(U_1,U_2)$ is of the type in the condition.
By the hypothesis, one may find a pseudoconvex open set  $\Omega$
such that
 $\widehat{W}\subset\Omega$ and $      U_2\cap\partial\Omega\not=\varnothing.$
 Let $U_3$ be the connected component of $U_2\cap \Omega$ which contains $U_1.$
 It is clear that  $\partial U_3\cap  U_2\cap\partial\Omega\not=\varnothing.$
 Therefore, by Cartan--Thullen Theorem, there are a  function
$f\in\mathcal{O}(\Omega),$   a sequence of points
$((z_j,w_j))_{j=1}^{\infty}\subset U_3$ and  a point $
(z_0,w_0)\in   \partial U_3 \cap  U_2\cap\partial\Omega$  with
$\vert f(z_j,w_j)\vert \to\infty$ and $(z_j,w_j) \to (z_0,w_0)$ as
$j\to\infty.$  Then $f$ restricted to $W$ extends to a holomorphic
function $\hat{f}$ on $\widehat{W}^{\text{o}}$ and  $f=\hat{f}$ on
$\widehat{W}^{\text{o}}.$
 By the above assumption there is $h\in\mathcal{O}(U_2)$ such that $h=\hat{f}$ on $U_1.$ This implies that
$f=h$ on $U_3.$ But  the latter identity contradicts the fact that
\begin{equation*}
\lim\limits_{j\to\infty}\vert f(z_j,w_j)\vert
=\infty\quad\text{and}\quad \lim\limits_{j\to\infty}\vert h
(z_j,w_j)\vert=\vert h(z_0,w_0)\vert.
\end{equation*}
 Hence, the proof is finished.
 \end{proof}
%\begin{lem}\label{lem2}
%   Let  $D\subset \C$ be an open set and  G,\ A, \  B $   Then
% $\widehat{W}^{\text{o}}$ is the envelope of holomorphy of $W.$
%  \end{lem}

 \smallskip

 Before we present the proof of Theorem A, we have to introduce some notation and terminology
 (see also Subsections 2.1--2.3 of  \cite{pn2})
% For $z\in \C^{N}$ and $r>0,$ let $\B(z,r)$ denote the ball with center $z$ and radius $r.$

 For two sets   $T\subset S,$  the {\it characteristic function} $1_{T,S}:\ S\longrightarrow\{0,1\}$ is given by:
 $1_{T,S}=1$ on $T$ and $1_{T,S}=0$ on $S\setminus T.$
 %denotes  the characteristic function of $T$ with respect to $S.$}

 Let us come back to the beginning of Subsection 1.3. For a curve $\mathcal{C}:=\{\gamma(t):\ t\in [0,1]\}$
 (resp. a closed curve which is the boundary of a Jordan domain $\Gamma:\ E\longrightarrow\C$),
 $\gamma:\ [0,1]\longrightarrow \C$ (resp. $\Gamma|_{\partial E}$) is said to be a {\it parametrization.}
 Moreover, two curves with corresponding parameterizations $\gamma_1,$  $\gamma_2$ are said to {\it have the
 same end-points} if  $\gamma_1(0)=\gamma_2(0)$ and  $\gamma_1(1)=\gamma_2(1).$
 %Now let us come back to the

 \smallskip

\noindent{\bf Proof of Theorem A.}
%  Now assume that  $A$ and $B$
%  do not possess any extendible points. We like to show that
%  $\widehat{W}^{\text{o}}$ is  the envelope of holomorphy of $W.$
 Since $D$     is Jordan-curve-like on the  set
  $A$ of positive length,   we may find a sequence $(A_k)_{k=1}^{\infty}$   of relatively
   open subsets of $\partial D$
     such that  $D$ is Jordan-curve-like      on $A_k,$ $k\geq1,$
  \begin{equation}\label{eq1_thmA}
   A\subset A_k,\ A_k\searrow A_0\ \text{as}\ k\to\infty,\qquad
   \text{and}\
 A_0\setminus A\ \text{is of zero length.}
   \end{equation}
   Using (\ref{eq1_thmA}) and  applying  Theorem 4.6 in \cite{pn2} yields  that
   \begin{equation} \label{eq2_thmA}
\omega(z,A,D)= \omega(z,A_0,D),\qquad z\in D   .
\end{equation}
Let $\mathcal{P}_D$ be the Poisson projection of $D$  (see
\cite[Subsection 4.3]{ra}). By  Theorem 4.3.3 in \cite{ra}, we
have
   \begin{equation*}
 \omega(\cdot,A_k,D) =\mathcal{P}_D[1_{\partial D\setminus A_k,\partial D}],\qquad  k\geq 0.
\end{equation*}
On the other hand, using (\ref{eq1_thmA}) and applying the
monotone convergence theorem yields that
 \begin{equation*}
 \lim\limits_{k\to\infty} \mathcal{P}_D[1_{\partial D\setminus A_k,\partial D}] =\mathcal{P}_D[1_{\partial D\setminus A_0,\partial D}].
\end{equation*}
This, combined with  (\ref{eq2_thmA}), implies that
 \begin{equation}\label{eq3_thmA}
\omega(z, A_k,D)\nearrow  \omega( z, A,D)\ \text{as}\ k
\to\infty,\ z\in D .
\end{equation}

Now we construct two sequences of open sets
$(D_k)_{k=1}^{\infty}\subset\C$ and
$(G_k)_{k=1}^{\infty}\subset\C$ with $D\subset D_k$ and $G\subset
G_k.$ To do this, for every $k\geq 1,$ write
\begin{equation}\label{eq4_thmA}
A_k=\bigcup\limits_{l\in I_k} A_{kl} ,\qquad
B_k=\bigcup\limits_{m\in J_k} B_{km}
\end{equation}
where  the $A_{kl}$   (resp. $B_{km}$)  are pairwisely disjoint
open connected subsets of $\partial D$ (resp. $\partial G$) and
the index set $I_k$  (resp. $J_k$) contains at most countably many
points.  Now we use the hypothesis that $D$  (resp. $G$) is
Jordan-curve-like at all points of $\overline{A}$  (resp.
$\overline{B}$). After a possible change of $A_k$   (resp. $B_k$)
we may assume
%Consequently, by removing at most countably many
%points from $\overline{A}_k$  (resp. $\overline{B}_k$) %and by
%taking  a conformal transformation (if necessary),
 that for every  fixed $k\geq 1,$ $(\overline{A_{kl}})_{ l\in I_k}$   (resp.
$(\overline{B_{km}})_{ m\in J_k}$) are pairwisely disjoint, they
are either curves or closed curves and their types are either $1$
or $2 .$ Moreover, using ( if necessary) a conformal
transformation: $z\mapsto \frac{1}{z-z_0}$ where $z_0$ is an
arbitrary point in $D,$ we may assume, without loss of generality
that $\partial D$ and $\partial G$ are bounded (hence, compact)
sets.
%\footnote{ An arc
%(resp. a closed curve) is, by definition, the image of a one-to-one continuous  map $\gamma:\
%[0,1]\longrightarrow\C$  (resp.  $\gamma:\ S^1 \longrightarrow\C,$
%where  $S^1:=\{z\in\C:\ \vert z\vert=1\}$).
% $\gamma$ is then called a parameterization of the arc (resp. closed curve).
% For an arc  with parametarization $\gamma,$
% the  open arc is the image $\gamma( (0,1))$
%and its two end-points are $\gamma(0)$ and $\gamma(1).$
%}

For   every  open (Jordan) curve
$A_{kl}$ (resp. $B_{km}$) of type $1$ % whose points  are of type 1  (see Subsection 2.1 in \cite{pn2} for the notion of types),
we fix a parametrization $a_{kl}:\
[0,1]\longrightarrow\overline{A_{kl}}$ (resp.   $b_{km}:\
[0,1]\longrightarrow\overline{B_{km}}$)  and
   find,
  using   a geometric argument which is  based on the compactness and the connectedness of $\overline{A_{kl}}$
  (resp.  $\overline{B_{km}}$),   an open  curve $ \Gamma_{kl}$  (resp.  $ \Lambda_{km}$)
    with a parametrization $\gamma_{kl}:\
  [0,1]\longrightarrow\overline{\Gamma_{kl}}$  (resp.  $\lambda_{km}:\
  [0,1]\longrightarrow\overline{\Lambda_{km}}$)   satisfying the following
  properties:
  \begin{itemize}
  \item[(a)] $\Gamma_{kl}$ has the same end-points as those of $A_{kl}$
  and  $\max\limits_{t\in [0,1]} \vert a_{kl}(t)-\gamma_{kl}(t)\vert<\frac{1}{k}$
  (resp.   $\Lambda_{km}$ has the same end-points as those of $B_{km}$
  and  $\max\limits_{t\in [0,1]} \vert b_{km}(t)-\lambda_{km}(t)\vert<\frac{1}{k}$);
  \item[(b)] $\overline{A_{kl}} \cup \overline{\Gamma_{kl}}$ is the boundary of a Jordan domain $\Delta_{kl}$
  (resp.   $\overline{B_{km}} \cup \overline{\Lambda_{km}}$ is the boundary of a Jordan domain $\Phi_{km}$);
    \item[(c)]  $\Delta_{kl}\cap D=\varnothing$  (resp.   $\Phi_{km}\cap
    G=\varnothing$).
      \end{itemize}

 For   every  closed (Jordan) curve $\overline{A_{kp}}$ (resp. $\overline{B_{kr}}$) of type $1$ % whose points  are of type 1  (see Subsection 2.1 in \cite{pn2} for the notion of types),
we fix a parametrization $a_{kp}:\ \partial E\longrightarrow
A_{kp}$ (resp.   $b_{kr}:\ \partial E\longrightarrow B_{kr}$)  and
   find,
  using as above   a geometric argument,   a  closed curve  $ \Gamma_{kp}$  (resp.  $ \Lambda_{kr}$)
    with a parametrization $\gamma_{kp}:\
  \partial E\longrightarrow \Gamma_{kp}$  (resp.  $\lambda_{kr}:\
  \partial E\longrightarrow \Lambda_{kr}$)   satisfying the following
  properties:
  \begin{itemize}
  \item[(d)] $A_{kp}\cap \Gamma_{kp}=\varnothing$ and
  $  \max\limits_{t\in \partial E} \vert a_{kp}(t)-\gamma_{kp}(t)\vert<\frac{1}{k}$
  (resp.    $B_{kr}\cap \Lambda_{kr}=\varnothing$
  and  $\max\limits_{t\in \partial E} \vert b_{kr}(t)-\lambda_{kr}(t)\vert<\frac{1}{k}$).
  Moreover,   $A_{kp} \cup \Gamma_{kp}$ is the boundary of an open annulus  $\Delta_{kp}$
  (resp.   $B_{kr} \cup\Lambda_{kr}$ is the boundary of an open annulus  $\Phi_{kr}$)
    such that $\Delta_{kp}\cap D=\varnothing$  (resp.   $\Phi_{kr}\cap G=\varnothing$).
     \end{itemize}

For   every   curve or closed curve $\overline{A_{kl}}$ (resp.
$\overline{B_{km}}$) of type $2$ let
\begin{itemize}
  \item[(e)] $\Gamma_{kl}:=\varnothing,$ $\Delta_{kl}:=\varnothing$ (resp.    $\Lambda_{km}:=\varnothing,$
  $\Phi_{km}:=\varnothing$).
     \end{itemize}

 Now we are able to  define the two sequences of open sets $(D_k)_{k=1}^{\infty}\subset\C$
and $(G_k)_{k=1}^{\infty}\subset\C.$ Indeed, for every $k\geq 1$
let
%consider the open set  $ D_k \subset\C$ (resp.  $G_k\subset \C$)   given by
   \begin{equation}\label{eq5_thmA}
   \begin{split}
  D_k &:= D\cup \Big ( \bigcup\limits_{l\in I_k^{'}}( \Delta_{kl}\cup A_{kl}) \Big)\cup
  \Big( \bigcup\limits_{p
  \in I_k^{''}}
   (\Delta_{kp}\cup A_{kp})\Big)  ,\\
 G_k &:=G\cup \Big ( \bigcup\limits_{m\in J_k^{'}}( \Phi_{km}\cup B_{km}) \Big)\cup
  \Big(\bigcup\limits_{r\in J_k^{''}}
   (\Phi_{kr}\cup B_{kr})\Big). %,\\
   %\widetilde{A}_k&:=A_k\cup \bigcup\limits_{l\in I_k}
   %\Gamma_{kl},\\
%\widetilde{B}_k&:=B_k\cup \bigcup\limits_{l\in J_k}
 %  \Lambda_{km},\\
   \end{split}
\end{equation}
where %$I_k$ (resp.  $J_k$) is as above, and
 $I^{'}_k$ (resp.
$J^{'}_{k}$)
    is the set of  all $l\in I_k$ (resp.  $m\in J_k$)  such that   $A_{kl}$ (resp.  $B_{km}$)  are open curves,  %whose points  are of type 1,
    and     $I^{''}_k$ (resp.   $J^{''}_{k}$)  is the set of  all $p\in I_k$ (resp.  $r\in J_k$)
           such that    $\overline{A_{kp}}$  (resp. $\overline{B_{kr}}$)  are closed curves.

  As a consequence of  the construction in (a)--(e) and  (\ref{eq4_thmA})--(\ref{eq5_thmA})
  we  are in position to show that
  \begin{equation}\label{eq6_thmA}
  \begin{split}
\omega(z,  A_k , D_k )&=
 \omega(z,A_k,D) %&=\omega(z,  \widetilde{A}_k , D_k ),
 \quad z\in   D;\\
\omega(w,  B_k , G_k )&=
 \omega(w,B_k,G),%&=\omega(w, \widetilde{B}_k , G_k )
  \quad  w\in G,\ k\geq 1.
\end{split}
\end{equation}
We only need to prove the first identity in (\ref{eq6_thmA}), the
other one can be proved similarly. In fact,  it suffices to show
that
\begin{equation}\label{eq6_supp_thmA}
\omega(z,  A_k , D)\leq \omega(z,A_k,D_k),\qquad z\in D,
\end{equation}
 since the converse inequality is
evident as $D\subset D_k,$ %$A_k\subset \widetilde{A}_k,$
$A_k\subset D_k$ (see (\ref{eq5_thmA})). To this end we define the
function $u:\ D_k\longrightarrow [0,1]$ as
\begin{equation*}
 u(z):=
\begin{cases}
\omega(z,A_k,D),
  & z\in   D,\\
  0, &z\in  D_k\setminus D.
  % A_k,\\
% \omega(z,A_{kl},\Delta_{kl}), & z \in \Delta_{kl},\ l\in I_k.\\
\end{cases}
\end{equation*}
Using \cite{ra}  we see that
 \begin{equation*}
\lim\limits_{z\to \zeta}\omega(z,  A_k , D )=0,\quad \zeta\in A_k. %;\qquad
%\lim\limits_{z\to z_0}\omega(z,  A_{kl} , \Delta_{kl} )=0,\quad z_0\in A_{kl},\
%l\in I_k.
\end{equation*}
Consequently, $u$ is subharmonic on $D_k,$ and  $u\leq
\omega(z,A_k,D_k),$ which implies (\ref{eq6_supp_thmA}).

Next, we like to check the hypothesis of Lemma \ref{lem1} in the
present context.  Therefore,
 fix a  nonempty open connected  set $U \subset \C^2$  such that
$U\not\subset \widehat{W}^{\text{o}}$  and $U\cap
\widehat{W}^{\text{o}}\not=\varnothing.$  We have to show that
   there is a pseudoconvex open set $\Omega$ in $\C^2$  such that
 $\widehat{W}\subset\Omega$ and $      U\cap\partial\Omega\not=\varnothing.$
In fact, we will choose  $\Omega$ as either
\begin{equation}\label{eq7_thmA}
\Omega=\Omega_k:=\widehat{\X}^{\text{o}}(A_k,B_k;D_k,G_k)\quad\text{or}
\quad  \Omega=D_k\times G_k
\end{equation}
for some $k\geq 1.$ In virtue of  (\ref{eq3_thmA}),
(\ref{eq5_thmA}),  and (\ref{eq6_thmA})--(\ref{eq7_thmA}), we
obtain that $\Omega_k$ is pseudoconvex and $\widehat{W}\subset
\Omega_k.$ Therefore, we only have to check that
 $      U\cap\partial\Omega\not=\varnothing.$
 Several cases are  to be considered.

 \noindent  {\bf Case I:}  $\left(U\cap\partial  \widehat{W}^{\text{o}}\right)\cap( D\times G)\not=\varnothing.$

 Fix a point $(z_0,w_0)\in  \left(U\cap\partial  \widehat{W}^{\text{o}}\right)\cap( D\times G).$
 Using the continuity of the harmonic measure, we see that
 \begin{equation*}
\partial  \widehat{W}^{\text{o}}\cap (D\times G)=\left\lbrace (z,w)\in D\times G:\ \omega(z,  A , D )+\omega(w,B,G)=1
\right\rbrace.
 \end{equation*}
In particular,  $ \omega(z_0,  A , D )+\omega(w_0,B,G)=1.$

Next, we fix  two points $(z_1,w_1)\in U\cap
\widehat{W}^{\text{o}}\cap( D\times G) $ and $(z_2,w_2)\in (U\cap(
D\times G))\setminus    \widehat{W}^{\text{o}} $  close to
$(z_0,w_0)$ such that
\begin{equation*}
 1< \omega(z_2,  A , D )+\omega(w_2,B,G),
  \end{equation*}
and
\begin{equation*}
\gamma:=\left\lbrace t(z_1,w_1)+(1-t) (z_2,w_2):\ t\in [0,1]
\right\rbrace \subset  U   \cap( D\times G).
 \end{equation*}
 Hence,
 \begin{equation*}
  \omega(z_1,  A , D )+\omega(w_1,B,G)<1< \omega(z_2,  A , D )+\omega(w_2,B,G)
  \end{equation*}
In virtue of (\ref{eq3_thmA}) and (\ref{eq6_thmA}), the
monotonically decreasing  sequence $ \left(\omega(\cdot,  A_k ,
D_k ) \right)_{k=1}^{\infty}$  (resp.   $ \left(\omega(\cdot,  B_k
, G_k ) \right)_{k=1}^{\infty}$)  of continuous functions
converges uniformly to $\omega(\cdot,A,D)$  (res.
 $\omega(\cdot,A,D)$) on
 some open neighborhoods of $z_1$ and $z_2$  (resp.  $w_1$ and $w_2$). Consequently, we may
 choose   $\Omega=\Omega_k$ for a sufficiently big
 $k$  in (\ref{eq7_thmA})
 such that
  \begin{equation*}
  \omega(z_1,  A_k , D_k )+\omega(w_1,B_k,G_k)<1< \omega(z_2,  A_k , D_k )+\omega(w_2,B_k,G_k)
  \end{equation*}
 This  implies  that  there is $(z_3,w_3)\in\gamma$ such that
    $ \omega(z_3,  A_k , D_k )+\omega(w_3,B_k,G_k)=1.$ Hence,
    $        U\cap\partial\Omega\not=\varnothing.$   So case I is complete.

 \noindent  {\bf Case II:}  $\left(U\cap\partial  \widehat{W}^{\text{o}}\right)\cap ( D\times G)=\varnothing.$

First we show that
\begin{equation}\label{eq1_Case_II}
 \partial  \widehat{W}^{\text{o}}\cap\Big ( (\partial D\setminus \overline{A})\times
 (\partial G\setminus  \overline{B})\Big)=\varnothing.
 \end{equation}
 Indeed, suppose  the contrary in order to get a contradiction
 and
 let  $(z_0,w_0)$ be an arbitrary  point in
the left hand side of (\ref{eq1_Case_II}).
  %$ (\partial D\setminus \overline{A})\times
 %(\partial G\setminus  \overline{B}).$
 %In virtue of    (\ref{eq1_thmA}) and
 %   the construction (a)--(d) above and the continuity of the harmonic measure,
 Since  $(z_0,w_0) \in (\partial D\setminus \overline{A})\times
 (\partial G\setminus  \overline{B}),$ we have  (see \cite{ra})
 \begin{equation*}
 \lim\limits_{z\to z_0} \omega(z, A,D)=1,\qquad  \lim\limits_{w\to w_0} \omega(w, B,G)=1,
  \end{equation*}
  which proves that  $(z_0,w_0) \not\in \overline{  \widehat{W}^{\text{o}}}. $
  Hence, we obtain the desired contradiction, and the proof of
  (\ref{eq1_Case_II}) is complete.

  %for sufficiently big $k.$ Fix such a $k$ and define $\Omega$ as in  in (\ref{eq7_thmA}).
  %Consequently,
  %there is a sufficiently small open neighborhood $U_3$ of $(z_0,w_0)$ such that $U_3\subset U_2$ and
 %\begin{equation*}
 %U_3\cap\partial  \widehat{W}^{\text{o}}=U_3\cap\partial  \Omega.
 %\end{equation*} Since
 %$U_2\cap\partial  \widehat{W}^{\text{o}}\not=\varnothing,$ it follows that $U_2\cap\partial  \Omega_k\not=\varnothing
 %.$ Hence, this subcase is finished.

Using (\ref{eq1_Case_II}), the obvious inclusion
$\left(U\cap\partial  \widehat{W}^{\text{o}}\right)\subset
\overline{D}\times \overline{G}$ and the assumption of  Case II,
we see that there are two subcases to consider.

{\bf Subcase}  $\left(U\cap\partial
\widehat{W}^{\text{o}}\right)\cap (\overline{A}\times
\overline{G}) \not=\varnothing.$ Let  $(z_0,w_0)$ be a point in
this intersection. Since  $D$ is Jordan-curve-like at all points
of $\overline{A},$  we may choose  a sufficiently large $k_0$ such
that
\begin{equation}\label{eq0_subcase}
\left\lbrace (z,w_0):\ \vert z-z_0\vert  %+\vert w-w_0\vert
<\frac{4}{k_0}\right \rbrace\subset U,
\end{equation}
and that all points of $T:=\left\lbrace z\in\partial D:\ \vert
z-z_0\vert<\frac{4}{k_0}\right \rbrace$ are either of the same
type 1 or of the same type 2.
 There are two subsubcases to consider.

   {\bf Subsubcase} The points of $T$ are  %or $B_{k_0m}$
     of type 1.   Since $z_0\in\overline{A},$   there exist $z_1\in A$ and $l\in I_{k_0}$  such that
  $\vert z_0-z_1\vert <\frac{1}{k_0}$ and
   $z_1\in A_{k_0l}.$
    Using (a)  above we can choose $t_1\in
(0,1)$ such that $z_1=a_{k_0l}(t_1).$ Now setting
$z_2:=\Gamma_{k_0l}(t_1),$  we see that
 $z_2\in\Gamma_{k_0l}$ and $\vert z_2-z_1\vert <\frac{1}{k_0}.$
 This, coupled with  (\ref{eq0_subcase}), gives that $(z_2,w_0)\in U.$

  Now we choose  $\Omega:=D_k\times G_k$ in (\ref{eq7_thmA}). It remains to show that
$(z_2,w_0)\in \partial\Omega.$ Since   $\Gamma_{k_0l}$
    is of type 1, it follows from (\ref{eq5_thmA}) that
$\Gamma_{k_0l}\subset \partial D_{k_0}.$ Hence, $z_2\in \partial
D_{k_0}.$     Therefore,  $(z_2,w_0)\in \subset
\partial\Omega.$

In summary, we have shown that $(z_2,w_0)\in U\cap\partial
\Omega.$
%If $z_0$ and $w_0$ are   of type 2, then clearly
%there is an open neighborhood $U_3\Subset U_2$ of $(z_0,w_0)$ such that
%$U_3\cap\partial  \widehat{W}^{\text{o}}=U_3\cap\partial  \Omega_k,$  $k\geq 1.$
Hence,   this subsubcase is completed.

 {\bf Subsubcase}  The points of $T$ are % and $B_{k_0m}$
  of type 2.
% Observe that the possibility that $z_0$ is an interior point of $A$ and $w_0$ is an interior point of $B$
% cannot happen. Otherwise,
%  \begin{equation*}
%    \lim\limits_{z\to z_0} \omega(z, A,D)= \lim\limits_{w\to w_0} \omega(w, B,G)=0,
%  \end{equation*}
%which shows that $(z_0,w_0)\in \widehat{W}^{\text{o}}.$

Since $z_0\in\overline{A},$   there exists $z_1\in A$  such that
  $\vert z_0-z_1\vert <\frac{1}{k_0}.$  For every $k\geq 1,$ let
  $l_k$ be the unique index in $I_k$  such that $z_1\in A_{kl_k}.$
  Since $(A_k)_{k=1}^{\infty}$ is a decreasing sequence of
  relatively open subsets
 of $\partial D,$  so is the sequence $(A_{kl_k})_{k=1}^{\infty}.$
 Put $H:=\bigcap\limits_{k=1}^{\infty} A_{kl_k}.$ We like to show that
 \begin{equation}\label{eq_thmA_subsubcase2}
H=\{z_1\}.
\end{equation}
Indeed, suppose in order to reach a contradiction that
$H\not=\{z_1\}.$  Then the interior of $H$ (in the relative
topology of $\partial D$) contains an open  (Jordan) curve
$\mathcal{C}$. On the other hand, we know from (\ref{eq1_thmA})
that $ A_k\searrow A_0$
 and $A_0\setminus A$ is of zero length, and it is easy to see
 that $H\subset A_0.$ Consequently,  $\mathcal{C}\setminus A$ is
 of zero length. Hence, by Definition 2  all points of $A\cap \mathcal{C}$ are
 extendible points of $A,$ this  contradicts the hypotheses of Theorem A. Hence,
(\ref{eq_thmA_subsubcase2})  has been proved.

In virtue of (\ref{eq_thmA_subsubcase2}) we may find a
sufficiently large $k$ such that $\sup\limits_{x,y\in
A_{kl_k}}\vert x-y\vert <\frac{1}{k_0}$ and $A_{kl_k}$ is an open
(Jordan) curve. Let $z_2$ be  an end-point  of
$\overline{A}_{kl_k}.$ So by the choice of $k,$ we have $\vert
z_1-z_2\vert <\frac{2}{k_0}.$ This, coupled with the previous
estimate $\vert z_0-z_1\vert <\frac{1}{k_0}$ and
(\ref{eq0_subcase}), gives that $(z_2,w_0)\in U.$

On the other hand, since $z_2$ is an end-point of $A_{kl_k},$ it
follows from (\ref{eq5_thmA}) that $z_2\in\partial D_k.$
 Now we choose
$\Omega:=D_k\times G_k$ in (\ref{eq7_thmA}). Consequently,
$(z_2,w_0)\in\partial \Omega.$

Summarizing, we obtain $(z_2,w_0)\in U\cap
\partial\Omega_{k}.$
 Hence, this subsubcase is completed.

{\bf  Subcase}   $\left(U\cap\partial
\widehat{W}^{\text{o}}\right)\cap (\overline{D}\times
 \overline{B})
\not=\varnothing.$ It is similar to the previous subcase.

Hence, the proof of the theorem is complete.
 \hfill $\square$

\smallskip

\noindent {\bf Proof of Theorem B.}
%For every $k\geq 1,$ let
%\begin{equation*}
%A_k:=D_k\setminus \overline{D},\qquad B_k:=G_k\setminus
%\overline{G}.
%\end{equation*}
Using this and hypothesis ($\Hc$)  and applying  Proposition 3.7
in \cite{pn3},  we see that, for every $k\geq 1,$
 \begin{equation}\label{eq0_thmB}
\begin{split}
\omega(\zeta,A,D_k)&= \lim\limits_{z\to\zeta,\  z\in D_k}
\omega(z,  A , D_k )=\lim\limits_{z\to\zeta,\  z\in D} \omega(z,A,D)=0,\qquad \zeta\in A,\\
  \omega(\eta,B,G_k)&=            \lim\limits_{w\to\eta,\ w\in G_k}\omega(w,  B , G_k )
  =\lim\limits_{w\to\eta,\ w\in G} \omega(w,B,G)=0,\qquad
                 \eta\in B.
\end{split}
\end{equation}
%and
%\begin{equation}\label{eq00_thmB}
%\omega(\zeta,A_k,D_k)=0,\   \omega(\eta,B_k,G_k)=0,      \qquad
%\zeta\in A_k,\eta\in B_k.
%\end{equation}
Consequently, arguing as in the proof (\ref{eq6_thmA}) one can
show that
\begin{equation}\label{eq1_thmB}
\omega(z,  A , D )= \omega(z,A,D_k),\quad z\in D,\qquad \omega(w,
B , G)= \omega(w,B,G_k),      \quad w\in G.
\end{equation}
Now we are able  to check the hypothesis of Lemma \ref{lem1} in
the present context.  To this end,
 fix a nonempty connected open set $U \subset  \C^{n+m}$ such that
$U\not\subset \widehat{W}^{\text{o}}$  and
$U\cap\widehat{W}^{\text{o}}\not=\varnothing.$ We will choose the
pseudoconvex open set $\Omega$ as
\begin{equation}\label{eq2_thmB}
\text{either}\quad
\Omega=\Omega_k:=\widehat{\X}^{\text{o}}(A,B;D_k,G_k)\quad\text{or}\quad
\Omega=D_k\times G_k,
\end{equation}
and need to verify that $\widehat{W}\subset \Omega$ and  $
U\cap\partial\Omega\not=\varnothing.$ Observe that  $\Omega$ is
pseudoconvex which follows from the definition of
$\widehat{\X}^{\text{o}}(A,B;D_k,G_k)$  and the fact that
 $D_k,$  $G_k$  are pseudoconvex.
On the other hand, using (\ref{eq0_thmB})--(\ref{eq2_thmB}) and
the fact that $A\subset D_k,$  $B\subset G_k,$ $D\subset D_k,$
$G\subset G_k,$ we easily see that $\widehat{W}\subset \Omega.$
Therefore, it remains to show that  $
U\cap\partial\Omega\not=\varnothing.$
% Since  $U\cap\partial
%\widehat{W}^{\text{o}}\not=\varnothing,$ it suffices to show that
%\begin{equation*}
%U\cap\partial \widehat{W}^{\text{o}}\subset U\cap\partial\Omega.
%\end{equation*}
To do this let $(z_0,w_0)$ be an arbitrary point in $
U\cap\partial \widehat{W}^{\text{o}}.$ There are several cases to
consider.

\smallskip

\noindent {\bf  Case I:} $(z_0,w_0)\in D\times G.$

\smallskip

Then there is a sequence  $((z_j,w_j))_{j=1}^{\infty}\subset
\widehat{W}^{\text{o}}$ such that
\begin{equation*}
\lim\limits_{j\to\infty}(z_j,w_j)=(z_0,w_0),\quad
\omega(z_0,A,D)+\omega(w_0,B,G) \geq 1.
\end{equation*}
This, combined with (\ref{eq1_thmB})--(\ref{eq2_thmB}), implies
that  $(z_0,w_0)\in\partial\Omega_k$ for arbitrary $k\geq1.$
Hence, choosing $\Omega:=\Omega_k$ for any $k\geq 1$ in
(\ref{eq2_thmB}) case I is completed.

\noindent {\bf  Case II:} $(z_0,w_0)\in \partial D\times
\overline{G}.$

Two subcases are to be considered.

\smallskip

\noindent {\bf  Subcase} $(z_0,w_0)\in A \times \overline{G}.$
%By
%(\ref{eq0_thmB}), we have
 %\begin{equation*}\lim\limits_{z\to z_0,\  z\in D}
 %\omega(z,A,D)=0.
%\end{equation*}
%Consequently, one may find  a sequence
%$((z_j,w_j))_{j=1}^{\infty}\subset (D\times G)\cap U$ such that
%\begin{equation*}
%\lim\limits_{j\to\infty}(z_j,w_j)=(z_0,w_0),\quad
%\lbrack\omega(z_j,A,D)+\omega(w_j,B,G)<1.
%\end{equation*}
%This, combined with (\ref{eq1_thmB})--(\ref{eq2_thmB}), implies
%that $(z_0,w_0)\in \overline{\Omega}_k.$
Recall from Definition 3 that $D_k\cap \overline{D}=D\cup A$ and
$\bigcap\limits_{k=1}^{\infty} D_k= D\cup A.$ Consequently, there
are an integer  $k$  and a point $z_1\in\partial D_k$ such that
$(z_1,w_0)\in U.$  Now put $\Omega:=D_k\times G_k.$ Then we  have
that $(z_1,w_0)\in U\cap\partial \Omega.$

\smallskip

\noindent {\bf  Subcase} $(z_0,w_0)\in (\partial D\setminus
A)\times \overline{G}.$ Since $z_0\in \partial D\setminus A,$ it
follows from Definition 3 that $z_0\in \partial D_k$ for $k\geq
1.$    Now choosing $\Omega:=D_1\times G_1,$  we see that
$(z_0,w_0)\in U\cap\partial \Omega.$

\smallskip

This completes case II.

\smallskip

\noindent {\bf  Case III:} $(z_0,w_0)\in \overline{ D}\times
\partial G.$  It is similar to case II.

\smallskip

 Hence, the proof of Theorem  B is finished. \hfill
$\square$

\smallskip

\noindent {\bf Proof of Proposition C.} First we consider the case
$\vert J\vert =1.$ Using  (i) and (iii) we may find a open
neighborhood $U$ of $\overline{A}$ in $\C^n$ which is relatively
compact   and a  $\mathcal{C}^2$ smooth strictly plurisubharmonic
defining function $\rho$ on $U$ such that
 $D\cap U=\left\lbrace z\in U:\ \rho(z)<0 \right\rbrace.$
  For every $z\in U\cap \partial D,$ let $v_z$ be the outward normal vector $v_z$ of $D$ at  $z.$
Using the smoothness in (iii) one may find
 a sufficiently small number $t_0>0$ such that the map
$\Theta:\  (U\cap \partial D)\times [0,t_0)\longrightarrow \C^n,$
given by
\begin{equation*}
\Theta(z,t):=z+ v_z,
\end{equation*}
is diffeomorphic onto the set $V\subset\C^n.$ Geometrically, $V$
is a tube with the base $U\cap \partial D$ and  with the height
$t_0.$

Since $\Theta\Big(((U\cap \partial D)\setminus A)\times
[0,t_0)\Big)$ is relatively closed in $V,$ there is a smooth
function $\lambda$ defined on $V$ such that $0\leq \lambda\leq 1$
on $V$ and
$$\{z\in V: \lambda(z)=0\}=\Theta\Big(((U\cap
\partial D)\setminus A)\times [0,t_0)\Big).$$
%and % $(D\cap U)\cup A= \overline{D}\cap \left\lbrace z\in U:\ \lambda(z)<0 \right\rbrace.  $
 %$ A= \partial D\cap \left\lbrace z\in U:\ \lambda(z)<0 \right\rbrace.  $
%Now one chooses the sequence of $\mathcal{C}^2$ smooth  strictly plurisubharmonic
%defining functions $(\rho_k)_{k=1}^{\infty}$ on $U$  given by
%\begin{equation*}
%\rho_k(z):=\rho(z)+\frac{1}{k}\cdot\lambda(z),\qquad z\in U.
%\end{equation*}
For every $k\geq 1,$ we define an open set
\begin{equation*}
D_{k}:=D\cup  \left\lbrace z\in U:\
\rho(z)-\frac{1}{k}\lambda(z)<0\right\rbrace.
\end{equation*}
Then using the above properties of $\rho,$ $\Theta$ and $\lambda,$
it can be checked that there exists  a sufficiently large  $N>0$
such that  $(D_{N+k})_{k=1}^{\infty}$ satisfies Definition 3.
Hence,
 the pair $(A,D)$   satisfies hypothesis ($\Hc$).

 The general case  (i.e. $\vert J\vert$ is at most countable)
may be done in  the same way  using (ii) and the fact that an
increasing union of pseudoconvex open sets is again pseudoconvex.
\hfill $\square$

\end{document}